\documentclass[b5paper, 12pt]{article}
\usepackage{amsfonts}
\usepackage{amsmath}
\usepackage{epsfig}
\usepackage{graphicx}
\usepackage{amssymb}

 \newcommand{\R}{\mathbb R}
 \newcommand{\C}{\mathbb C}
 \newcommand{\N}{\mathbb N }
 \newcommand{\Si}{\mathrm{Si}}
 
\begin{document}

\begin{center}
\fontsize{14pt}{16pt} \bf On regularizations of the integral representation of Dirac delta function -- elementary approach
\end{center}

\begin{center}
\sf{ Grzegorz M. Koczan $^{1)}$, Piotr Stachura $^{2)}$  }
\end{center}

\medskip

{\small \centerline{ $^{1)}$ Department of Mechanical Processing of Wood}, 
     \centerline{ Warsaw University of Life Sciences-SGGW}
  \centerline{ ul. Nowoursynowska 159, 02-776 Warszawa, Poland}
 \centerline{ {\tt grzegorz\_koczan1@sggw.edu.pl}}
 \medskip
 
\centerline{ $^{2)}$ Department of Applied Mathematics}, 
     \centerline{ Warsaw University of Life Sciences-SGGW}
  \centerline{ul. Nowoursynowska 159, 02-776 Warszawa, Poland}
\centerline{ {\tt piotr\_stachura1@sggw.edu.pl}}
\bigskip

\begin{quote} {\bf Abstract.}
  {\it Dirac's Delta and its most popular integral representation are peculiar objects that often cause
    controversy over their mathematical status. The problem with handling these objects is a fact and does not contradict
    the existence of the well established theory of distributions. However, this situation may indicate the lack of sufficient
    popularization of the theory or its difficulties, or it may suggest the existence of hidden and poorly
    understood aspects of the theory. The authors therefore decided to present, in an elementary way, but with mathematical
    precision and without harm to the intuition, the path from the integral representation to the Dirac delta,
    starting with Schwartz's functional approach. Next, the considered representation is
    presented in a more intuitive sequential approach, formulated by Mikusi\'nski and Sikorski.
    Finally, we present the third regularization that can be related to Sato's approach to distributions.
  }
\end{quote}

\bigskip
\normalsize
\section{Introduction}

In many physical and technical applications, the following representation of the Dirac delta is used:
\small
\begin{equation}\label{delta:def1}
  \delta_{x_0}(x):=\delta(x-x_0):=\frac{1}{2 \pi} \int_{-\infty}^{\infty}\mathrm{d k\,} e^{ik(x-x_0)} 
\end{equation}
\normalsize
The above form of the formula is sometimes attributed to Augustin Louis Cauchy~(1789--1857), while its equivalent form is due to
Jean Baptiste Joseph Fourier~(1768--1830):
\small
\begin{equation}\label{delta:def2}
  \delta(x-x_0):=\frac{1}{2 \pi} \int_{-\infty}^{\infty}\mathrm{d k\,} \cos(k(x-x_0))
\end{equation}
\normalsize

The expression~(\ref{delta:def1}) was popularized by the English physicist Paul \\
Dirac~(1902--1984)~\cite{dirac1}.  The general theory of similar quantities was first clarified by the
French mathematician Laurent Schwartz~(1915--2002) in the
language of functionals~\cite{schw1, schw2, schw3} and a little later by Polish mathematicians Jan Mikusi\'nski~(1913--1987) and
Roman Sikorski~(1920--1983) in a language of function sequences~\cite{mik48, mik55, mik-sik57, mik-sik61, sik54}.
The above formulas also have a deeper justification as inverse Fourier transforms  of harmonic function or, in particular,
the constant function $1$. In this context, the sense of an integral expressions (\ref{delta:def1}) or  (\ref{delta:def2})
should not be categorically challenged. However, it should be remembered that the above integral has a
generalized sense, and its result is a distribution that goes beyond the set of functions.
Nevertheless, the question arises how to treat rigorously such an expression and how to elementarily
prove that it has the assumed properties.

In the historical context, it is often recognized that the first general and strict theory of distribution (the functional one)
was formulated  by L. Schwartz in the years 1945--1951. However, already in the mid-1930s, the Russian mathematician Sergei Sobolev~(1908--1989)
considered the concept of the weak derivative~\cite{sob36, sob38, kuta1}, which is crucial for the functional approach.
Mikusi\'nski and Sikorski, on the other hand, in the years 1948--1961 proposed an alternative approach to the theory with sequences.
However, it is much more difficult to decide the priority of the discovery of distributions as such.
One can get the impression that historians of science are slightly polarized here along national lines English-French or
physics and mathematics.  In this way, as a pioneer can be credited Paul Dirac, who introduced in 1926~\cite{dirac2}
the useful concept of delta, now called by his name. Dirac interpreted his delta function as a generalization of the Kronecker
delta for continuous sets.

Over a hundred years earlier, in 1822, in a work on the equation of thermal conductivity
and its solutions~\cite{fou22}, Fourier already used formulas similar to these:
\small
\begin{align}\label{delta:fou1}
    f(x)&=\frac{1}{2 \pi} \int_{-\infty}^{\infty}\mathrm{d k\,} \int_{-\infty}^{\infty}\mathrm{d x_0\,}\cos(k(x-x_0)) f(x_0)\\
      \label{delta:fou2}
    f(x)&=\int_{-\infty}^{\infty}\mathrm{d x_0} \frac{\sin(k(x-x_0))}{\pi(x-x_0)} f(x_0)\,,\,k\rightarrow\infty
\end{align}
\normalsize
On the one hand, it can be said that Fourier did not write explicitly the formula~(\ref{delta:def2}),  but on the other hand, his advantage
remained mathematical accuracy and usage of well  defined concepts. So he  implicitly used the formula (\ref{delta:def2}) as
if inside the reverse (Fourier) transform (\ref{delta:fou1}). Under the sign of an additional integral,
the order of which he rearranged, he computed the indefinite integral (\ref{delta:def2}) and obtained (\ref{delta:fou2}).
However, he stopped there, because he could not calculate the  limit, which in the usual  sense did not exist.
In Fourier's calculations,  in expressions of type (\ref{delta:fou1}) or  (\ref{delta:fou2}) always  appeared transformed  function,
which, together with the additional integral allowed formally to compute the limit and return to the initial function
(in the sense of the truth of the equation (\ref{delta:fou2}), see[6]).

In this method, one can even see an
outline of  Schwartz's later theory. Similar conclusions to Fourier were reached in 1827 by Cauchy~\cite{cauchy}, who gave
the complex equivalent of the integral (\ref{delta:fou1}). According to the source~\cite{katz}, Cauchy had some other ideas about
singular functions as early as 1815. In addition  to this early contribution of French scholars,  there were  contributions of
English scholars George Green~(1793--1841) and Oliver Heaviside~(1850--1925).
The first of them considered functions of which Laplacian had the properties of the Dirac delta in the present sense; the other 
used the derivative of a step function (often called the Heaviside function), which is in fact a delta distribution.

This article is not a brief introduction to the theory of distributions, nor does it present all the underlying ideas.
It only shows, {\em using the most elementary methods possible}, how to interpret a certain heuristic
and useful formula in two approaches to  the theory of distributions.

\section{Regularization in terms of Schwartz's distributions}
   
\newcommand{\sD}{{\mathcal D}}
Let  ${\mathcal D}:=C_0^\infty(\R)$ be the (vector) space of smooth, complex functions with compact support.
Each function from $\sD$ has derivatives of arbitrary degree and is identically $0$ off some bounded interval $[\alpha,\beta]$,
depending on the function, an easy construction of an example is presented in the Appendix. Elements of ${\mathcal D}$ will be called  {\em test functions}.
This is a topological space, but we need not and will not define this topology here (see e.g. \cite{schw2}).
It is sufficient, for our purposes, to define what it means for a seqence of test functions to converge to zero.
Let $(f_n)$ be a sequence of test functions. According to Schwartz, such a sequence is converging to the zero function if the supports of 
all $f_n$'s are contained in {\em a fixed} interval $[-M,M]$ and the sequence, and all its derivatives converge uniformly:
\small
\begin{equation}\label{zbiezwD}   \frac{d^kf_n(x)}{d x^k}\rightrightarrows 0\,,\,k=0, 1, 2, \dots
\end{equation}
\normalsize
A double arrow in this article will mean {\em uniform convergence}~\cite{mik-sik57}.
The condition that supports  of all functions $f_n$ are contained  in the specified interval is important:
if $f\in \sD, f(1)\neq 0$ then the sequence $f_n(x):=\frac{1}{n}f(\frac{x}{n})$ {\em does not converge} in  $\sD$,
although it converges uniformly with all derivatives.

Let us also consider a broader set of locally integrable functions
$L^1_{loc}(\R)$, i.e. for any real numbers  $\alpha<\beta$ there exists the integral $\int_\alpha^\beta \mathrm{d x\,} \varphi(x)$.
Therefore, to each such function can be assigned a linear functional $\varphi[\,]$ on  $\sD$:
\small
\begin{equation}\label{funk:defi}
\varphi[f]:=\int_{-\infty}^\infty \mathrm{d x\,} \varphi(x) f(x)=\int_{-M}^M \mathrm{d x\,} \varphi(x) f(x),
\end{equation}
\normalsize
where $M$ is such that the support of  $f$ is contained  in $[-M, M]$.
In order not to complicate notation, the functional $\varphi[\,]$ is distinguished from the function $\varphi(\,)$
only by the shape of parentheses. This functional is continuous in the following sense: if the sequence of test functions $f_n$ tends to
zero according to~(\ref{zbiezwD}), then the (numeric) sequence $\varphi[f_n]$ tends to zero.

The idea of distributions  is to generalize functions to functionals. All linear and continuous functionals on
$\sD$ in the above sense are called distributions and denoted $\sD'$, and the functionals of type (\ref{funk:defi})
form its proper subspace, because
it turns out that for a  distribution (functional) $T\in\sD'$ there may not be a function $\varphi(\,)\in L^1_{loc}(\R)$,
for which $T=\varphi[\,]$. The best example is the Dirac delta described by the generalized  function (\ref{delta:def1}): its corresponding
functional should  meet the following condition:
\small
\begin{equation}\label{delta-defi}
  \delta[f]=f(0)   
\end{equation}
\normalsize
The linear and continuous functional $\delta[\,]$ is specified correctly, but  its ``integral kernel'', i.e. the default (generalized)
function $\delta(\,)$ is no longer an ordinary object. Nevertheless, this study shows how to treat this object in the
context of the representation (\ref{delta:def1}). The most literal understanding of such a
generalized function is made possible by the Mikusi\'nski-Sikorski sequential  approach. However, we will start
with a description oriented towards the functional approach of  (\ref{delta:def1}).

In the first step, the improper integral will be rewritten as the limit (for $x_0=0$):
\small
\begin{equation}\label{delta-R}
  \delta(x):=\lim_{R\rightarrow\infty}\frac{1}{2\pi} \int_{-R}^{R}\mathrm{d k\,} e^{ikx}
\end{equation}
\normalsize
It turns out that the adoption of symmetrical boundaries here is not crucial for further proceedings, but
it significantly shortens the notation. Rewriting (\ref{delta:def1}) in the form of the  limit (\ref{delta-R}) does not change much,
because such a limit in the ordinary sense still does not exist.
So let's switch the order of
taking the limit  $R\rightarrow\infty$ and transition from function to functional:
\small
\begin{equation}\label{delta-R-f}
  \delta_R(x):=\frac{1}{2\pi} \int_{-R}^{R}\mathrm{d k\,} e^{ikx}\,\rightarrow\delta_R[\,]
  \stackrel{R\rightarrow\infty}{\longrightarrow} \delta[\,]
\end{equation}
\normalsize
In this procedure, with the changed order of taking the limit, all objects (function, functional and limit functional)
have a well-defined mathematical sense. In fact, the description using the family of functions is very close to the sequential approach
(Fig. 1, Fig. 3), but  while the limit of the family of  functionals $\delta[\,]$ is not a problem, understanding the limit
for the family of functions $\delta_R(\,)$ is rather problematic (Fig.~\ref{fig:2}).
\begin{figure}[h]
\begin{center}
\includegraphics[scale=0.7]{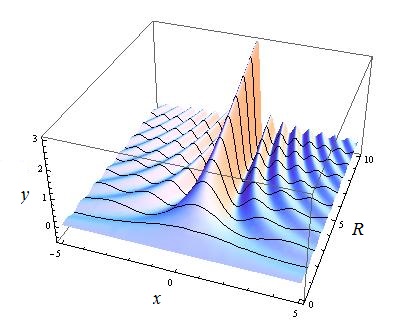}
\caption{\small Three-dimensional visualization of the $\delta_R(x)$ family of functions describing the Dirac delta with marked contour
  lines corresponding to the integer  values of the parameter $R>0$. }\vspace{-3ex}
\label{fig:1}
\end{center}
\end{figure}

However, it turns out, as shown by  Mikusi\'nski and Sikorski (see the next section), that this  is possible by
considering the limit of the family (sequence)
of primitive functions to $\delta_R(\,)$. The sense of this limit will be clarified and its existence will be
proven in the second section of the article. In this  section we treat  the formula (\ref{delta:def1}) according to the
procedure (\ref{delta-R-f}) and prove that it satisfies the functional definition (\ref{delta-defi}) of the Dirac delta.

Let us start with the elementary definite integral:
\small
\begin{equation}\label{delta-R-sin}
  \delta_R(x)=\frac{1}{2\pi} \int_{-R}^{R}\mathrm{d k\,} e^{ikx}=\frac{\sin(R x)}{\pi x}\,,x\neq 0,\,\,\delta_R(0)=\frac{R}{\pi}
\end{equation}
\normalsize
\begin{figure}[h]
\begin{center}
\includegraphics[scale=0.6]{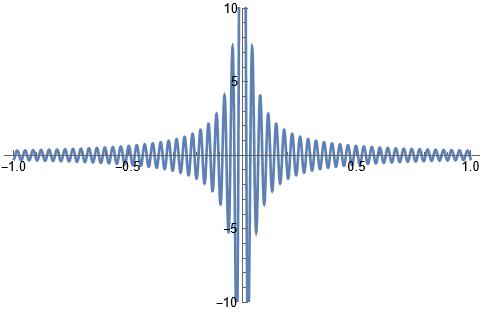}
\caption{\small Graph of one of the functions describing the Dirac delta for a conventionally large value  of the parameter
    $R=n=180$. Values of such functions oscillate in the range $[-(\pi x)^{-1}, (\pi x)^{-1}]$ and do not go to zero.}\vspace{-3ex}
\label{fig:2} 
\end{center}
\end{figure}
\noindent
and consider the functional corresponding to the function (\ref{delta-R-sin}):
\small
\begin{equation}
   \delta_R[f]=\int_{-M}^M\mathrm{d x\,} \frac{\sin(R x)}{\pi x} f(x)
\end{equation}
\normalsize

We would like to compute the limit of this expression as $R\rightarrow\infty$. Let us begin by rewriting  the expression as:
\small
\begin{equation}\label{delta-R-f1}
  \begin{split}
  \delta_R[f]&=\int_{-M}^M\mathrm{d x\,} \frac{\sin(R x)}{\pi x} f(x) =\\
 &=  \frac{1}{\pi}\int_{-M}^M\mathrm{d x\,} \sin(R x)\frac{ f(x)-f(0)}{x}+\frac{f(0)}{\pi}\int_{-M}^M\mathrm{d x\,} \frac{\sin(R x)}{x}
  \end{split}
\end{equation}
\normalsize
To calculate this limit, we will use several lemmas. The first of them requires only the application of l'Hospital's rule:

\noindent
{\bf Lemma 1.}\label{lemat1} {\it If $f\in\sD$ is a test function, then the function:
  $$g(x):=\frac{f(x)-f(0)}{x}\,,x\neq 0\,,\,g(0):=\lim_{x\rightarrow 0} g(x)=f'(0)$$ has a continuous first derivative.
  (In fact, it is even smooth, but we will not need this fact; note that  unless $f(0)=0$, its support is not compact)}\\
\hspace*{\fill}$\Box\quad\quad$

\noindent
{\bf Lemma 2.} {\it Let a function $g$ be continuously differentiable, then:
\small $$ \lim_{R\rightarrow\infty}\int_{-M}^M \mathrm{d x\,} g(x) \sin(Rx) = 0$$}\vspace{1ex}

\noindent
{\em Proof:\,} Integration by parts. \hspace*{\fill}$\Box\quad\quad$

\noindent Using the above lemmas, we obtain the limit of~(\ref{delta-R-f1}):
\small
$$\lim_{R\rightarrow \infty}\delta_R[f]=\frac{f(0)}{\pi} \lim_{R\rightarrow\infty}\int_{-M}^M\mathrm{d x\,} \frac{\sin(R x)}{x}=
\frac{f(0)}{\pi} \int_{-\infty}^\infty\mathrm{d x\,} \frac{\sin(x)}{x}$$
\normalsize
It remains to calculate the improper integral (convergent in the sense of the Riemann integral, divergent in the sense of
the Lebesgue integral). An elementary but tedious calculation of this so-called Dirichlet integral can be found e.g. in the 
Fichtenholz textbook~\cite{ficht2}, it can also be calculated using the integration with respect to a parameter or the residuum method.
Continuing  the elementary nature of this work, we will use Fubini's Theorem in its simplest version.

\noindent
{\bf Lemma 3.} {\it 
  The following improper Riemann integral is convergent to $\pi$.
 \small \begin{equation}\label{calka-n}
     \int_{-\infty}^\infty\mathrm{d x\,} \frac{\sin(x)}{x}=\pi
  \end{equation}
}\vspace{1ex}

\noindent
{\em Proof:} It is sufficient to show  that  {\small $\displaystyle  \int_0^\infty\mathrm{d x\,} \frac{\sin(x)}{x}=\pi/2$}.
We integrate the function $h(x,\alpha):=e^{-\alpha x}\sin(x)$ on the square $K:=[0,R]\times [0,R]$; let
$\displaystyle S(R):=\int_K h(x,\alpha)$. Integrating by parts twice, we obtain:
\small
\begin{equation}\label{SR}\begin{split}
  S(R)&=\int_0^R \mathrm{d \alpha\,} \int_0^R \mathrm{d x\,} e^{-\alpha x}\sin(x)=\\&=
  \int_0^R \frac{ \mathrm{d \alpha}}{1+\alpha^2}- \int_0^R \mathrm{d \alpha\,} e^{-\alpha R}\frac{\cos R+\alpha \sin(R)}{1+\alpha^2}
  \end{split}
\end{equation}
\normalsize
Since  $\frac{|\alpha|}{1+\alpha^2}\leq \frac12$, the function in the second integral can be estimated by:
\small
$$ \left|  e^{-\alpha R}\frac{\cos R+\alpha \sin(R)}{1+\alpha^2}\right|\leq e^{-\alpha R}(1+\frac12),$$
\normalsize
and we get:
\small
$$ \left|\int_0^R \mathrm{d \alpha\,} e^{-\alpha R}\frac{\cos R+\alpha \sin(R)}{1+\alpha^2}\right|\leq
\frac32 \int_0^R \mathrm{d \alpha\,} e^{-\alpha R}=\frac{3(1-e^{-R^2})}{2 R}$$
\normalsize
therefore {\small $\displaystyle \lim_{R\rightarrow\infty}S(R)=\lim_{R\rightarrow\infty}\arctan(R)=\pi/2$. }
On the other hand:
\small
\begin{equation*}
S(R)=\int_0^R \mathrm{d x\,} \int_0^R \mathrm{d \alpha\,} e^{-\alpha x}\sin(x)=
  \int_0^R \mathrm{dx\,}\frac{\sin(x)}{x}- \int_0^R \mathrm{d x\,} \frac{e^{-R x}\sin(x)}{x}
\end{equation*}
\normalsize
Since  $|\frac{\sin(x)}{x}|\leq 1$ the second integral tends to zero (as above) and we get:
\small
$$\pi/2=\lim_{R\rightarrow\infty}S(R)=\lim_{R\rightarrow\infty}\int_0^R \mathrm{d x\,}\frac{\sin(x)}{x}$$
\normalsize\hspace*{\fill}$\Box\quad\quad$

And, finally, for any test function $f\in\sD$:
\small
\begin{equation}\label{delta-ost} 
  \lim_{R\rightarrow\infty}\delta_R[f]=f(0)=\delta[f]
\end{equation}
\normalsize                                      
This way, the defining property of the Dirac delta function (\ref{delta-defi}) has been proven in the Schwartz functional approach.

\section{Regularization in the sequential approach of Mi\-ku\-si\'ns\-ki and  Sikorski}

So far, it has been shown that the integral formula (\ref{delta:def1}) should, in principle, be treated as a family of
integrals (functions) (\ref{delta-R-sin}). This family is not convergent  in the space of ordinary functions, however,
the family of functionals  corresponding to it according to the recipe (\ref{funk:defi}) is convergent.
Schematically, this has been  presented in (\ref{delta-R-f}), where  ``a delay'' in passing to the limit is explicitly mentioned:
firstly, one replaces a function by a corresponding functional and then passes to the limit.
Such a "delayed" limit should be treated not as an error, but as an element of the definition
of a distribution (as a functional). This delay, equivalent to the change of the order of integration, is the key aspect of Schwartz's
theory, and the need for it was already recognized by Fourier in (\ref{delta:fou1}) and (\ref{delta:fou2}).
However, let us ask ourselves whether the Dirac delta can be treated not as a functional,
but as a generalized  limit that could be expressed by the formula (\ref{delta:def1})?
Let us now sketch the approach, the idea is similar to the construction of real numbers from rationals by Cauchy sequences.

Let $\varphi_n:\R\rightarrow\C$ be continuous functions. The sequence of functions $(\varphi_n)$ will be called
{\em a fundamental sequence}, if there is a sequence of functions $\Phi_n$ and an integer $k\geq 0$ such that:

\small
\begin{equation}\label{ciag-fun-defi}
  \begin{split}
    \Phi_n^{(k)}(x)&:=\frac{d^k\Phi_n(x)}{d x^k}=\varphi_n(x)\,{\rm \, and}\\
    \Phi_n(x)& \rightrightarrows \Phi(x)\,\,\mbox{{\rm \,on \, any\,  interval \,}} [a,b]
  \end{split}
\end{equation}
\normalsize
The limit function $\Phi(x)$ is then continuous. A sequence of functions converging uniformly on each interval $[a,b]$
will be called convergent {\em almost uniformly}. It follows directly from the definition that if the
sequence $(\varphi_n)$ is fundamental and derivatives of the order $l$, $\varphi_n^{(l)}$ are continuous,
then the sequence $(\varphi_n^{(l)})$ is fundamental as well.

Two fundamental sequences are {\em equivalent} iff sequences of their primitives (of order $k$) are almost
uniformly convergent to the same function.
\small
\begin{equation}\label{ciagi-row}
  \begin{split} (\varphi_n) \sim (\psi_n) \iff  \Big( & \exists\,k\in \N\cup \{0\}: \varphi_n=\Phi_n^{(k)}\,,\,\psi_n=\Psi_n^{(k)}
    \mbox{{\rm \,and\,}}  \\
  &  \Phi_n\rightrightarrows \Phi \leftleftarrows \Psi_n\, \mbox{{\rm \,on \, any\,  interval \,}}  [a,b] \Big)
  \end{split}
\end{equation}
\normalsize
This is really an equivalence relation and its equivalence classes {\em are distributions} in the sequential approach of
Mikusi\'nski and Sikorski; the equivalence class of a sequence $(\varphi_n(x))$ is denoted by $\langle\varphi_n(x)\rangle$. 

If the sequence of continuous functions converges almost uniformly, its class is the same as the class of the constant sequence
consisting of the limit function. In this sense, every continuous function represents a distribution, much like the Cauchy
sequence of rational numbers converging to a rational number is equivalent  to a constant sequence composed of its limit.

Let's give an example of equivalent  sequences: the sequence $\Phi_n(x):=-\frac{\cos(nx)}{n}$ is almost
uniformly convergent to $0$, Thus, it represents the distribution $\langle 0\rangle$; let $\Psi_n(x):=0$.
We then have $\Phi_n\rightrightarrows 0\leftleftarrows\Psi_n$. In addition, the sequence $\varphi_n(x)=\Phi_n''(x)=n\cos(nx)$
is fundamental, and $(n\cos(nx))\sim (\Psi_n''(x))$ thus  $\langle n\cos(nx)\rangle=\langle 0\rangle$. 

Let's do a calculation which shows the relationship of fundamental sequences with the functional approach. Let a  sequence
$\varphi_n=\Phi_n^{(k)}$ be a fundamental sequence, and $\Phi_n$ converges almost uniformly to $\phi$.
For a test  function $f$, by performing successive integrations by parts,  we obtain: 
\small
\begin{equation}
  \begin{split}
    \varphi_n[f]&=\int_{-M}^{M} \mathrm{d x \,} \Phi_n^{(k)}(x) f(x)=- \int_{-M}^{M} \mathrm{d x \,} \Phi_n^{(k-1)}(x) f'(x)=\\
    &=\int_{-M}^{M} \mathrm{d x \,} \Phi_n^{(k-2)}(x) f''(x)=\dots=(-1)^k \int_{-M}^{M} \mathrm{d x \,} \Phi_n(x) f^{(k)}(x)
  \end{split}
\end{equation}
\normalsize
therefore  \small
$\displaystyle  \lim_{n\rightarrow\infty}\varphi_n[f]=(-1)^k \int_{-M}^{M} \mathrm{d x \,} \Phi(x) f^{(k)}(x)=(-1)^k\Phi[f^{(k)}]$.
\normalsize
By the definition of convergence in  $\sD$, the last expression defines a continuous functional of $f$ i.e. a distribution.

Going back to $\delta(x)$, we want to show that the family of continuous functions (\ref{delta-R-sin}) for natural $R=n$ (fig.~\ref{fig:3}),
i.e. 
\small
$$\delta_n(x)=\frac{1}{2\pi} \int_{-n}^{n}\mathrm{d k\,} e^{ikx}=\frac{\sin(n x)}{\pi x}\,,x\neq 0,\,\,\delta_n(0)=\frac{n}{\pi},$$
\normalsize
is a fundamental sequence that has the desired properties of the Dirac delta.
\begin{figure}[h]
\begin{center}
\includegraphics[scale=0.7]{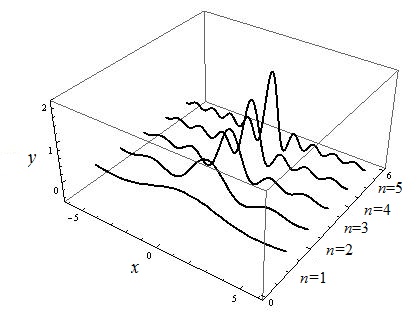}
\caption{\small Graphs of the first five functions of the fundamental sequence $\delta_n$ describing the Dirac delta.
  All functions of this sequence are smooth, even though continuity is sufficient
  in the fundamental sequence.}\vspace{-3ex}
\label{fig:3} 
\end{center}
\end{figure}
Therefore, it is necessary to find a sequence of primitive functions that converge almost uniformly. The first integration leads
to the integral sine, written with a capital letter:
\small
\begin{equation}\label{teta-n}
  \theta_n(x):=\int_0^x\mathrm{d t\,} \delta_n(t)=\int_0^x\mathrm{d t\,} \frac{\sin(nt)}{\pi t}=
  \frac{1}{\pi}\int_0^{nx} \mathrm{d t\,} \frac{\sin(t)}{t}=: \frac{1}{\pi}\Si(nx)
\end{equation}
\normalsize
\begin{figure}[h]
\begin{center}
\includegraphics[scale=0.7]{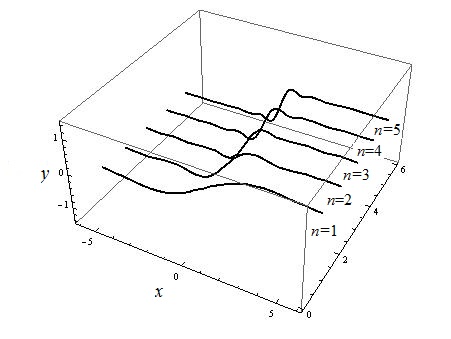}
\caption{\small Graphs of the first five functions $\theta_n$  of the fundamental sequence approximating  the step function.}\vspace{-3ex}
\label{fig:4} 
\end{center}
\end{figure}
Compared to $\delta_n$, $\theta_n$ has at least a point limit which, under (\ref{calka-n}), is a step function (Fig.~\ref{fig:4}):
\small
\begin{equation}\label{teta}
  \lim_{n\rightarrow\infty}\theta_n(x)=\left\{\begin{array}{rr}-\frac12 & x<0\\0 & x=0\\ \frac12 & x>0\end{array}\right\}=:\theta(x)
\end{equation}
\normalsize
\begin{figure}[h]
\begin{center}
\includegraphics[scale=0.7]{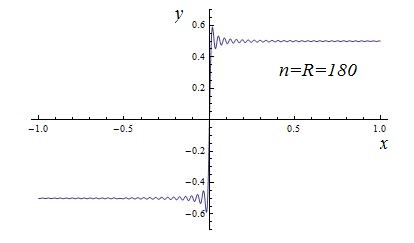}
\caption{\small Graph of one of the functions describing the step function  for a conventionally large value  of the parameter
  $R=n=180$. Taking a closer look at the flat fragments of the graph, one can see that the  derivatives do not converge
  in the usual sense to $0$  as $n$ increases (these derivatives are  of course the functions $\delta_n(x)$).}
\label{fig:6}
\end{center}
\end{figure}
Therefore, the sequence $\theta_n$ does not converge almost uniformly and we need second integration (by parts):
\small
\begin{equation}\label{Delta}
  \Delta_n(x):=\frac{1}{\pi} \int_0^x\mathrm{d t\,}\Si(nt)=\frac{x}{\pi}\Si(nx)+\frac{\cos(nx)-1}{n \pi}
\end{equation}
\normalsize
\begin{figure}[h]
\begin{center}
\includegraphics[scale=0.7]{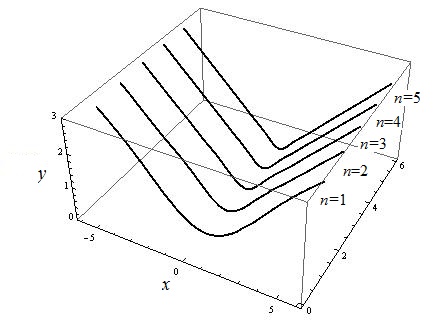}
\caption{\small Graphs of the first five functions  of the fundamental sequence $\Delta_n(x)$
  approximating the half of the absolute value function.}\vspace{-3ex}
\label{fig:5}
\end{center}
\end{figure}
We calculate the point limit in the same way as before (\ref{teta}) and this time we get a continuous function (fig.~\ref{fig:5}):
\small
\begin{equation}\label{deltan}
  \lim_{n\rightarrow\infty}\Delta_n(x)=\frac{|x|}{2}=:\Delta(x) 
\end{equation}
\normalsize
Let us put the thesis that this  is a uniform convergence in the following:

\noindent
{\bf Lemma 4.} {\it The sequence $\delta_n(x):=\frac{sin(nx)}{\pi x}$ is a fundamental sequence for which the sequence of second-order
  primitive functions $\Delta_n(x)$ given by (\ref{Delta}) tends uniformly to the function $\Delta(x):=\frac{|x|}{2}$.}\vspace{1ex}

\noindent
{\em Proof:\,} In order to prove uniform convergence, it is useful to express the integral sine in a slightly different way:
\small
\begin{equation}
  \Si(nx)=\int_0^{nx} \mathrm{d u\,}\frac{\sin(u)}{u}=\int_0^{\frac{nx}{2}} \mathrm{d y\,}\frac{\sin(2 y )}{y}=
  \int_0^{\frac{nx}{2}} \mathrm{d y\,}\frac{2 \sin(y )\cos (y)}{y}
\end{equation}
\normalsize
Integrating by parts, we obtain:
\small
\begin{equation}\label{Si-n}
  \Si(nx)=\frac{2 \sin^2(\frac{nx}{2})}{n x}+\int_0^{\frac{nx}{2}}\mathrm{d y\,} \frac{\sin^2(y)}{y^2}=\frac{1-\cos(nx)}{n x}+
  \int_0^{\frac{nx}{2}}\mathrm{d y\,} \frac{\sin^2(y)}{y^2}
\end{equation}
\normalsize
This way, by (\ref{teta-n}) and  (\ref{calka-n}), we get the value of the following improper integral:
\small
\begin{equation}\label{sin-2}
  \int_0^\infty\mathrm{d y\,}\frac{\sin^2(y)}{y^2}=\frac{\pi}{2}    
\end{equation}
\normalsize
Using (\ref{Si-n}), the sequence of second-order primitive functions $\Delta_n$ (\ref{Delta}) gets a more concise form:
\small
$$\Delta_n(x)=\frac{x}{\pi}\int_0^{\frac{nx}{2}}\mathrm{d y\,}\frac{\sin^2(y)}{y^2}$$
\normalsize
Since $\Delta_n$ and $\Delta$ are even, it is sufficient to estimate the difference of  $\Delta_n(x)$ and $\Delta(x)$ for $x\geq 0$;
given (\ref{sin-2}) we get:
\small
\begin{equation}
  \begin{split}
    \left|\Delta_n(x)-\Delta(x)\right|&=
    \frac{x}{\pi} \left|\int_0^{\frac{nx}{2}}\mathrm{d y\,}\frac{\sin^2(y)}{y^2}  -\frac{\pi}{2}\right| 
    =   \frac{x}{\pi} \int_{\frac{nx}{2}}^\infty\mathrm{d y\,}\frac{\sin^2(y)}{y^2} \leq\\
    &\leq \frac{x}{\pi} \int_{\frac{nx}{2}}^\infty\frac{ \mathrm{d y}}{y^2} =\frac{2}{n\pi},   
  \end{split}
\end{equation}
\normalsize
therefore,  the convergence is uniform and the sequence $\delta_n(x)$ is a fundamental sequence, so it represents a distribution.
\hspace*{\fill}$\Box\quad\quad$

All that remains is to show that this distribution is the  Dirac delta.
In the second  section, it was shown that the integration of terms of  the fundamental
sequence $\delta_n(x)$ with a test  function, after passing to the limit leads to the desired value~(\ref{delta-ost}).
Now it is appropriate to list other properties that also define the Dirac delta. Such a property can be considered
a derivative in the distributional sense of the step function (with a unit height step), which is a derivative of the function
$\frac{|x|}{2}$ (for $x\neq 0$).

Here we touch on the problem of defining the derivative of a distribution in a sequential  approach.
We will not define this notion in general, and  will limit ourselves to the specific case where our distribution is
represented by a fundamental sequence $(\varphi_n)$ of continuously differentiable  functions.
Then, as already mentioned, $(\varphi_n')$ is also a fundamental sequence, which suggests the following definition of the
derivative of a  distribution:
\small
\begin{equation}\label{dist-poch}
  \langle\varphi_n\rangle':=\langle\varphi_n'\rangle
\end{equation}
\normalsize
If $(\psi_n)$ is another representation with  $\psi_n'$ continuous, then looking at the definition of equivalence~(\ref{ciagi-row}),
we can immediately see that $(\varphi_n')\sim(\psi_n')$, so the definition~(\ref{dist-poch}) does not depend on the choice of a
differentiable representative. It turns out that there is always such a representative (even a smooth one),
we refer the insightful reader to work~\cite{mik-sik57}. In our situation, we have:
$$\delta_n(x)=\theta_n'(x)=\Delta_n''(x)$$
The sequence $\theta_n$ is fundamental because it is an integral of a fundamental sequence,
so it represents a distribution $\langle \theta_n(x)\rangle$, while the sequence $\Delta_n$ converges almost uniformly,
i.e. it represents the distribution $\langle \Delta(x)\rangle$. We can therefore write:
$$\delta(x):=\langle \delta_n(x)\rangle=\langle \theta_n(x)\rangle'=\langle \Delta(x)\rangle''$$
Of course, the symbol $\delta(x)$ reintroduced here is neither the uniform nor the point limit but represents a distribution.

In the functional approach, where there is an easy identification of locally integrable functions with corresponding distribution,
we can  easily verify that $\delta(x)=\theta'(x)$ (after defining derivative of distribution), but in the sequential approach we  can't
immediately write $\langle\delta_n(x)\rangle=\langle\theta(x)\rangle'$, since $\theta(x)$ is not continuous, and  the
symbol $\langle\theta(x)\rangle$ has no meaning. This is why we had to choose a longer path to $\Delta(x)$,
in fact, also in this approach it is also possible to identify distributions with locally integrable functions, see~\cite{mik-sik57}).

\section{Regularization related to Sato's approach}
So far, both in Schwartz's and Mikusi\'nski-Sikorski's approach, regularization has been carried out by means of limiting procedure in
the improper integral. Limit  was performed in functional or sequential-distribution approach.
Now we will approach the regularization of this integral a little differently. This is related to the third approach to the
theory of distributions developed in 1958--1960 by Japanese mathematician Miki Sato~(1928--2023). It is based on boundary values of
holomorphic functions, but since our presentation was promised to be  {\em elementary}, after this short remark,
we refer to \cite{sato} for more information.

Let us divide the integral~(\ref{delta:def1}) into positive and negative frequencies including the converging factor:
\small
\begin{equation}\label{delta-sato}
  \begin{split} \frac{1}{2 \pi} \int_{-\infty}^{\infty}\mathrm{d k\,} e^{ik x}&:=
    \frac{1}{2 \pi} \left( \int_{-\infty}^0\mathrm{d k\,} e^{ik x}+ \int_0^\infty \mathrm{d k\,} e^{ik x}\right):=\\
    &=\frac{1}{2 \pi} \lim_{\epsilon\rightarrow 0^+}\left( \int_{-\infty}^0\mathrm{d k\,} e^{ik x}e^{k\epsilon}+
      \int_0^\infty \mathrm{d k\,} e^{ik x}e^{-k\epsilon} \right), 
  \end{split}
\end{equation}
\normalsize
the last integrals converge for $\epsilon>0$. Let's introduce notation
\small 
$$\tilde{\delta}^{\epsilon}(x):=\frac{1}{2 \pi}\int_{-\infty}^\infty\mathrm{d k\,} e^{ik x}e^{-|k|\epsilon}
=\frac{1}{2 \pi} \left( \int_{-\infty}^0\mathrm{d k\,} e^{ik x}e^{k\epsilon}+
  \int_0^\infty \mathrm{d k\,} e^{ik x}e^{-k\epsilon} \right)$$ 
\normalsize
The last two  integrals are elementary, and performing them, we get:
$$\tilde{\delta}^\epsilon(x)=\frac{1}{2 \pi i}\left(\frac{1}{x-i\epsilon}-
  \frac{1}{x+i\epsilon}\right)=\frac{\epsilon}{\pi}\frac{1}{x^2+\epsilon^2}.$$
This way, we have the family of functionals $\tilde{\delta}^\epsilon[\,]$ and, for $n\in \N$, we define the sequence
$\tilde{\delta}_n(x):=\tilde{\delta}^\epsilon(x)\,,\,\epsilon:=\frac1n$  i.e.
\small
\begin{equation}\label{delta-tilde}
  \tilde{\delta}_n(x)=\frac{n}{\pi}\frac{1}{1+n^2 x^2}.
\end{equation}
\normalsize
\begin{figure}[h]
\begin{center}
\includegraphics[scale=0.4]{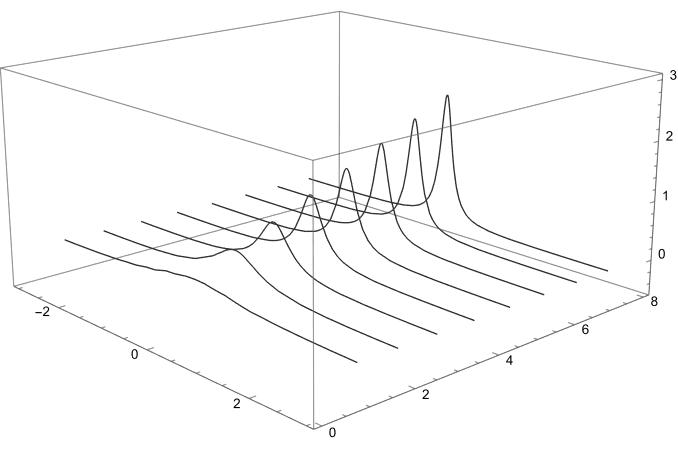}
\caption{\small Graph of several first functions $\tilde{\delta}_n(x)$. Their behaviour is  ``much more regular'' then that  of $\delta_n(x)$.}
\label{fig:7}
\end{center}
\end{figure}

\noindent
{\bf Lemma 5. } {\it The sequence  $(\tilde{\delta}_n)$ is fundamental and $\langle\tilde{\delta}_n\rangle=\langle\delta_n\rangle$.
  Moreover  for a function $f\in\sD$:
  $\quad\displaystyle \lim_{\epsilon\rightarrow 0^+}\tilde{\delta}^\epsilon[f]=f(0).$}

\noindent
{\em Proof: }Computation for $(\tilde{\delta}_n)$ is simpler than for $(\delta_n)$ since all integrals are elementary:
\small \begin{equation*}
  \begin{split}
    \tilde{\theta}_n(x)&:=\frac{1}{\pi}\int_0^x\mathrm{d t\,} \frac{n}{1+n^2 t^2}=\frac{1}{\pi}\arctan(nx),\\
    \tilde{\Delta}_n(x) &:=\frac{1}{\pi}\int_0^x \mathrm{d t}\,\arctan(nt)=\frac{x\arctan(n x)}{\pi}-\frac{\ln(1+n^2 x^2)}{2 \pi n}
    \stackrel{n\rightarrow\infty}{\longrightarrow} \frac{|x|}{2}
  \end{split}
\end{equation*}
\normalsize
We left to the reader to verify that the sequence  $(\tilde{\Delta}_n)$  converges almost uniformly. This way we get that the sequence
$(\tilde{\delta}_n)$ is fundamental and $\langle\tilde{\delta}_n\rangle=\langle\delta_n\rangle$.

To prove the second part, let us compute:
\small
\begin{equation*}
  \begin{split}
    \tilde{\delta}^\epsilon[f]&=\frac{1}{\pi}\int_{-M}^{M} \mathrm{d x\,}\frac{\epsilon}{x^2+\epsilon^2} f(x)=\\
    &=\frac{1}{\pi}\int_{-M}^{M} \mathrm{d x\,}\frac{\epsilon x}{x^2+\epsilon^2}\frac{f(x)-f(0)}{x}+
    \frac{f(0)}{\pi}\int_{-M}^{M} \mathrm{d x\,}\frac{\epsilon }{x^2+\epsilon^2}
  \end{split}
\end{equation*}
\normalsize
The second integral is elementary and we obtain $2\arctan(\frac{M}{\epsilon})$, and \\
\small $\displaystyle \lim_{\epsilon\rightarrow 0^+}2\arctan(\frac{M}{\epsilon})=\pi,$ \normalsize
so it remains to show that the limit of the first integral is $0$. Let $g(x):=\frac{f(x)-f(0)}{x}$, by Lemma 1,
this is a continuous function, and let $S$ be the supremum of $|g(x)|$ on $[-M,M]$.
The absolute value of the first integral is majorized  by:
\begin{equation}\label{limlog}
  \frac{2 S}{\pi}\int_{0}^{M} \mathrm{d x\,}\frac{\epsilon x}{x^2+\epsilon^2}=
  \frac{S\epsilon}{\pi}\left(\ln(M^2+\epsilon^2)-\ln(\epsilon^2) \right)\stackrel{\epsilon\rightarrow 0^+}{\longrightarrow} 0
  \end{equation}
\hspace*{\fill}$\Box\quad\quad$

The sequence $\tilde{\delta}_n(x)$ is well known as a representation of the Dirac delta, it is much more ``regular'' then $\delta_n(x)$, e.g.
it is converging to $0$ at all points $x\neq 0$ and its  primitives $\tilde{\theta}_n(x)$ and $\tilde{\Delta}_n(x)$ are elementary functions.
\section{What values  does the Dirac delta attain?}
It is usually said that $\delta(x)=0$ for $x\neq 0$. To verify  such an equality, it is necessary to discuss {\em the value of a distribution at a point}.
Such a definition was given by {\L}ojasiewicz in \cite{loj1957}. It uses the convergence of distributions, which  would lead us
beyond the scope of this article.
We will present other, more elementary, arguments to support this claim.

In both approaches, distributions are generalization of functions and one may identify some distributions with functions:
in the functional approach, the equality $T=\varphi$ for $T\in \sD'$ and $\varphi\in L^1_{loc}(\R)$ holds iff $T(f)=\varphi[f]$ for any $f\in \sD$; 
and in the sequential approach $\langle \varphi_n\rangle= \varphi$, for continuous $\varphi$, iff $\varphi_n$ is equivalent to a constant sequence $\varphi$
(here again one may identify $\langle \varphi_n\rangle$ with a locally integrable function, see \cite{mik-sik57}).

In both approaches, one can also restrict
distributions to open sets in $U\subset \R$: in the functional approach, just consider test functions with support in $U$, and in the sequential
approach consider fundamental sequences of functions defined on $U$. If distributions are equal on the larger set, they are also equal on
the smaller one. In this sense, it may happen that although a distribution is not equal to any function on $\R$ it can be equal to some on an open
subset $U\subset \R$. This is the case of the Dirac delta:\vspace{1ex}

\noindent
{\bf Lemma 6. } \textit{On the open set $\R\setminus\{0\}$:  $\delta(x)=0$.}\\
{\em Proof:} a) The functional approach: from the definition, if $f\in \sD$ is supported in $\R\setminus\{0\}$ then 
$\delta[f]=f(0)=0$.\\
b) The sequential approach: it is easier to use $\tilde{\delta}_n$ representative.
The claim will follow from   almost uniform convergence of $\tilde{\delta}_n$ to $0$ on $\R\setminus\{0\}$. But this is  clear since 
$\tilde{\delta}_n(x)$ converges pointwise  to $0$ on $\R\setminus\{0\}$ and for $|x|\geq a$ there is the inequality 
$|\tilde{\delta}_n(x)|\leq |\tilde{\delta}_n(a)|.$

Let us also notice that on $\R\setminus\{0\}$, we may write $\langle \theta_n\rangle= \langle \theta\rangle$,
this makes sense since on this set the step function $\theta$ is continuous (even smooth). In fact, on this set the sequence
$\theta_n\rightrightarrows \theta$; it is sufficient to consider $x>0$:
\small
\begin{equation*}
  \begin{split}
    \left|\theta_n(x)-\frac12\right|&=\frac1{\pi}\left|\int_0^{nx} \mathrm {d t\,} \frac{\sin(t)}{t}-\frac{\pi}{2}\right|=
    \frac1{\pi}\left|\int_0^{nx} \mathrm{d t\,} \frac{\sin(t)}{t}-\int_0^{\infty} \mathrm {d t\,} \frac{\sin(t)}{t}\right|=\\
    &= \frac1{\pi}\left|\int_{n x}^{\infty} \mathrm{d t\,} \frac{\sin(t)}{t}\right|
  \end{split}
\end{equation*}
\normalsize
Since the integral \small$\displaystyle \int_0^{\infty} \mathrm{d t\,} \frac{\sin(t)}{t}$ \normalsize is convergent, for any $\epsilon>0$ there exists
$A>0$ such that for $x>A$:  \small$\displaystyle \int_x^{\infty} \mathrm{d t\,} \frac{\sin(t)}{t}\leq \epsilon $\normalsize.
So given $a>0$ and $\epsilon>0$, for any $n\geq \frac{A}{a}$ one has \small $ \left|\theta_n(x)-\frac12\right|\leq\frac{\epsilon}{\pi}$ \normalsize 
for any $x\geq a$,  what proves almost uniform convergence of $\theta_n$ (on $\R\setminus\{0\}$). A reader may also check, it is even easier, that
on this set also the sequence $(\tilde{\theta}_n)$ converges almost uniformly.
\hspace*{\fill}$\Box\quad\quad$
\section{Summary}
It was shown how to understand the divergent  Fourier-Cauchy integrals~(\ref{delta:def1}, \ref{delta:def2}) as representations of the Dirac delta, either
as a limit of functionals $\delta_R[\,]$~(\ref{delta-R-f1}), in Schwartz's sense or as a fundamental sequence
$\delta_n(x)$ in Mikusi\'nski-Sikorski approach, or using yet another type of regularization~(\ref{delta-sato}),
as a limit of $\tilde{\delta}_n(x)$~(\ref{delta-tilde}) (in both approaches). We also indicated how to understand that $\delta(x)=0$ for $x\neq 0$.
Moreover, let us note that Dirac delta is used to describe the point charge in electrodynamics  and, relatively recently, was used to describe
(much more complicated) point mass in General Relativity \cite{kata}.
Finally, refering to the title of this volume,  as it was kindly pointed out by an  anonymous
Referee, most of integrals and pointwise limits appearing in the article, e.g. (\ref{calka-n}, \ref{sin-2}, \ref{deltan}, \ref{limlog}),
can be computed directly by Mathematica (ver. 14.2.1), however, our goal in this article was to increase understanding, not computational skills.

\section{Appendix: Compactly supported smooth function on $\R$}
In this appendix, we present an easy construction of a smooth function with  compact support. In fact, for any real numbers $\alpha<\beta<\gamma<\delta$,
we construct a smooth, non-negative function with support in $[\alpha,\delta]$ which is one on $[\beta,\gamma]$. Let us define\vspace{-1ex}
$$f(x):=\left\{\begin{array}{lcl} \exp(-1/x) & {\rm for} & x>0\\0 & {\rm for } & x\leq 0\end{array}\right.\vspace{-1ex}$$
By l'Hospital's rule (and some induction), one proves that $f$ is a smooth function and all its derivatives at $0$ are $0$.
Let $g(x):=f(-x)$, and for $\alpha,\beta$ let $f_\alpha(x):=f(x-\alpha)$ and $g_\beta(x)=g(x-\beta)$ (see Fig.\ref{fig:8}).
Notice that for $\alpha<\beta$ the function $f_\alpha+g_\beta$ is strictly positive.
\begin{figure}[h]
\begin{center}
\includegraphics[scale=0.3]{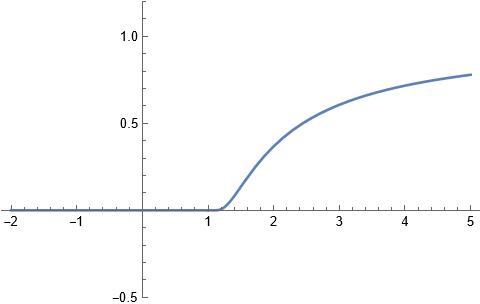}\,\qquad\,
\includegraphics[scale=0.3]{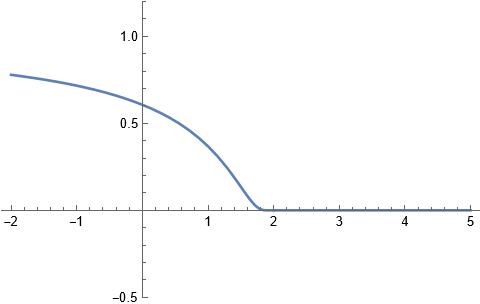}
\caption{\small Graphs of the functions $f_1(x)$ (left) and $g_2(x)$ (right).}
\label{fig:8}
\end{center}
\end{figure}\vspace{-2ex}

\noindent Define functions:
\small
$$\mathrm{F}_{\alpha,\beta}(x):=\frac{f_\alpha(x)}{f_\alpha(x)+g_\beta(x)}\quad,\qquad\mathrm{G}_{\gamma,\delta}(x):=\frac{g_\delta(x)}{f_\gamma(x)+g_\delta(x)}$$
\normalsize
These are smooth, unit ``up'' and ``down'' step functions on $[\alpha,\beta]$ and $[\gamma,\delta]$ intervals respectively, so their product 
$\mathrm{F}_{\alpha,\beta}(x)\mathrm{G}_{\gamma,\delta}(x)$ has the desired properties (see Fig.\ref{fig:9}). 
\begin{figure}[h]
\begin{center}
\includegraphics[scale=0.3]{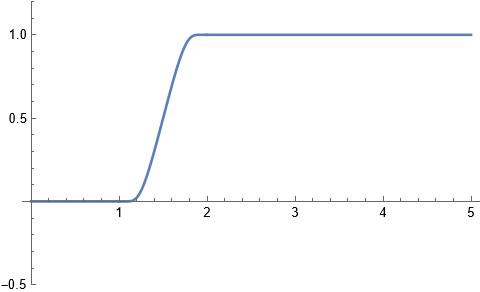}\,\,
\includegraphics[scale=0.3]{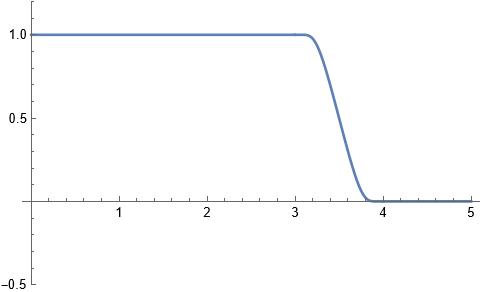}\,\,
\includegraphics[scale=0.3]{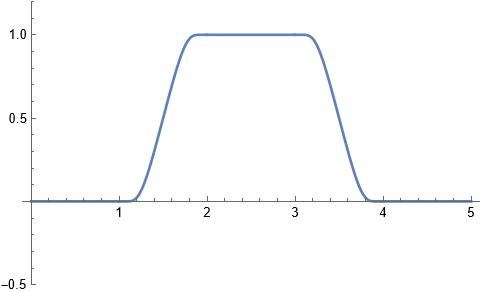}
\caption{\small Graphs of the functions $F_{1,2}(x)$ (left) and $G_{3,4}(x)$ (center) and their product (right).}
\label{fig:9}
\end{center}
\end{figure}
\pagebreak

\end{document}